\documentclass[10pt]{article}
\usepackage[tbtags]{amsmath}
\usepackage{epsfig,amstext,amssymb,amsthm,latexsym}
\allowdisplaybreaks[4]
\usepackage{amssymb,color}

\definecolor{c20}{rgb}{0.,0.7,0.}
\definecolor{c30}{rgb}{0.,0.,1.}
\definecolor{c40}{rgb}{1,0.1,0.7}
\definecolor{c50}{rgb}{1,0,0}

\def\ass#1{\textcolor{c50}{#1}}
\def\ass#1{#1}
\newtheorem{theo}{Theorem}[section]
\newtheorem{sat}[theo]{Proposition}
\newtheorem{de}[theo]{Definition}
\newtheorem{lem}[theo]{Lemma}
\newtheorem{exxa}[theo]{Example}

\newtheorem{korr}[theo]{Corollary}

\newtheorem{remarks}[theo]{Remarks}

\newcommand{\nelem}[1]{{Lemma \ref{#1}}}

\newcommand{\netheo}[1]{{Theorem \ref{#1}}}

\newcommand{\prooftheo}[1]{ \textsc{Proof of Theorem} \ref{#1} }

\newcommand{\prooflem}[1]{\textsc{Proof of Lemma} \ref{#1}}

\newcommand{\kb}[1]{\boldsymbol{#1}}
\newcommand{\vk}[1]{\kb{#1}}

\def\fracl#1#2{\biggr(\frac{#1}{#2} \biggl) }

\newcommand{\E}[1]{\mbox{\rm$\vk{E}$}\{#1\}}

\newcommand{\pk}[1]{\mbox{\rm$\vk{P}$} \{#1\} }

\newcommand{\pb}[1]{\mbox{\rm$\vk{P}$}\Bigl \{#1 \Bigr \}}

\newcommand{\R}{\!I\!\!R}

\newcommand{\inr}{\in \R}

\newcommand{\ldot}{,\ldots,}

\newcommand{\limit}[1]{\lim_{#1 \to   \infty}}

\newcommand{\BQN}{\begin{eqnarray}}
\newcommand{\EQN}{\end{eqnarray}}
\newcommand{\BQNY}{\begin{eqnarray*}}
\newcommand{\EQNY}{\end{eqnarray*}}

\newcommand{\BS}{\begin{sat}}
\newcommand{\ES}{\end{sat}}
\newcommand{\BT}{\begin{theo}}
\newcommand{\ET}{\end{theo}}
\newcommand{\BK}{\begin{korr}}
\newcommand{\EK}{\end{korr}}

\newcommand{\BD}{\begin{de}}
\newcommand{\ED}{\end{de}}

\newcommand{\BIT}{\begin{itemize}}
\newcommand{\EIT}{\end{itemize}}
\newcommand{\BDI}{\begin{description}}
\newcommand{\EDI}{\end{description}}

\newcommand{\BRM}{\begin{remarks}}
\newcommand{\ERM}{\end{remarks}}

\newcommand{\QED}{\hfill $\Box$}

\newcommand{\IF}{\infty}
\newcommand{\BTH}{\begin{theo}}
\newcommand{\ETH}{\end{theo}}
\newcommand{\BPR}{\begin{sat}}
\newcommand{\EPR}{\end{sat}}

\newcommand{\BEX}{\begin{exxa}}
\newcommand{\EEX}{\end{exxa}}

\newcommand{\BC}{\begin{cases}}
\newcommand{\EC}{\end{cases}}
\newcommand{\COM}[1]{}
\newcommand{\BL}{\begin{lem}}
\newcommand{\EL}{\end{lem}}

\def\SI{\Sigma}

\def\x{\vk{x}}
\def\y{\vk{y}}
\def\X{\vk{X}}

\def\U{\vk{U}}

\def\1d{\{1 \ldot d\}}
\def\njk{\{1 \ldot k\}}

\def\FRE{\mbox{Fr\'{e}chet }}

\newcommand{\equaldis}{\stackrel{d}{=}}

\begin{document}

\begin{center}
\thispagestyle{empty}

{\Large  On the residual dependence index }

{\Large of elliptical distributions }

       \vskip 0.4 cm
         \centerline{\large Enkelejd Hashorva}
\COM{        
\centerline{\textsl{Allianz Suisse }}
        \centerline{\textsl{CH-3001 Bern, Switzerland}}

}

\centerline{\textsl{Department of Mathematical  Statistics and Actuarial Science}}
        \centerline{\textsl{University of Bern, Sidlerstrasse 5}}
        \centerline{\textsl{CH-3012 Bern, Switzerland}}
       \centerline{\textsl{enkelejd.hashorva@stat.unibe.ch }}

\today{}

\end{center}

{\bf Abstract:} The residual dependence index of bivariate Gaussian
distributions is determined by the correlation coefficient. This tail index is of certain statistical importance when extremes and related rare events  of bivariate samples with asymptotic independent components are being modeled.
In this paper we calculate the partial residual dependence indices of a multivariate
elliptical random vector assuming that the associated random radius is in the Gumbel max-domain of attraction.
Furthermore, we discuss the estimation of these indices when the associated random radius possesses a Weibull-tail distribution.



\section{Introduction}
Let $(X_1,X_2)$ be a bivariate elliptical random vector with stochastic representation
\BQN \label{e} (X_1,X_2) &\equaldis &
R\Bigl(U_1, \rho U_1+ \sqrt{1-\rho^2} U_2 \Bigr), \quad \rho\in (-1,1),
\EQN
where the positive random radius $R$ is independent of $(U_1,U_2)$ which is uniformly distributed on the unit circle of $\R^2$. Here
$\equaldis$ stands for  equality of distribution functions.
A canonical example of a bivariate elliptical random vector is when $R^2$ is Chi-square  distributed, which implies
that $X_1,X_2$ are standard Gaussian random variables with mean 0, variance 1
 and correlation coefficient $\rho:=\E{X_1X_2}$.
It is well-known (see e.g., Reiss and Thomas (2007)) that in the
Gaussian model the correlation coefficient $\rho$ does not influence
the asymptotic dependence of the components. Roughly speaking this
means that the sample extremes of Gaussian random vectors are
asymptotically independent.

A tractable extension of the Gaussian model is the elliptical one, where $R$ is some general positive
random variable with distribution function $F$. By Lemma 12.1.2 of Berman (1992)
\BQN
\label{eq:berm12}
X_1 \equaldis X_2 \equaldis R U_1 \equaldis R U_2
\EQN
implying that distribution function of $X_1$ (denoted below by $Q$) is continuous.
Clearly, the joint dependence function of $X_1$ and $X_2$ is influenced by $\rho$. \\
Several authors have considered elliptical distributions for modelling of rare events with specific applications in insurance and finance.
Recent contributions in this directions are Peng (2008), Li and Peng (2009). In the first mentioned article the author deals with the situation that the marginal distributions are regularly varying implying that the components are asymptotically dependence. In such a model, also considered in Kl\"uppelberg et al.\ (2007), novel estimation techniques are presented.\\
 In Berman (1983,1992), Hashorva (2005a,b), Abdous et al.\ (2005,2008)  further probabilistic results are obtained when $F$ is in the Gumbel max-domain of attraction which implies that the components $X_1$ and $X_2$ are asymptotically independent. \\
 In this paper, with motivation from above mentioned contributions,
 we focus on the quantification of the asymptotic dependence of elliptical random vectors.

An interesting measure of the asymptotics dependence (see e.g., Peng (1998, 2007), de Haan and Ferreira (2006), Reiss and Thomas (2007)) is the function $\chi(u)$ defined by
\BQNY
\chi(u):= \frac{\pk{X_1>u,X_2>u}}{\pk{X_1>u}}\in [0,1], \quad u>0.
 \EQNY
If for some constant $c\in [0,1]$ we have
\BQN \label{eq:ind.a}
\lim_{u\to \IF } \chi(u)&=& c\in (0,1],
\EQN
then $X_1$ and $X_2$ are said to be asymptotically dependent. In our setup of bivariate elliptical random
vectors with stochastic representation \eqref{e} this is the case when
 $R$ has distribution function $F$ in the \FRE max-domain of attraction
(or equivalently, $F$ is regularly varying with positive index $\gamma$).
See Berman (1992) or Hashorva (2005a,2006b) for further details.
Important statistical applications can be found in Kl\"uppelberg et al.\ (2007).

In both other cases of max-domain of attraction, i.e., $F$ is in the Gumbel or the Weibull max-domain of
attraction we have (see Hashorva (2005a)) $c=0$,
which means that $X_1$ and $X_2$ are asymptotically independent.\\
In extreme value theory asymptotic independence is a nice property,
however, $c=0$ in \eqref{eq:ind.a} merely means that  $\pk{X_1>u,X_2>u}$ converges faster to 0 than the marginal survival probability $\pk{X_1>u}$ (if $u\to \IF$).\\
One successful approach to model the asymptotic independence is the estimation of the residual dependence index $\eta\in (0,1]$ (see e.g.,
Peng (1998, 2007, 2008), de Haan and Peng (1998), or  de Haan and Ferreira (2006)). We note in passing that recent ideas in testing asymptotic can be found in H\"usler and Li (2009). Now, information about $\eta$ is available if for any
$x,y$ positive
\BQN \label{defSU}
S_u(x,y):= \frac{
\widetilde{S}_u(x,y)}{ \widetilde{S}_u(1,1)} &\to & S(x,y)\in (0,\IF),
\quad u \to \IF,
  \EQN
with
$$   \widetilde{S}_u(x,y):= \pk{Q(X_1)> 1- x/u, Q(X_2)> 1- y/u}, \quad  u>0,$$
since for any $c>0$ and for some $\eta\in (0,1]$ we have the important scaling relation
\BQNY
 S(cx,cy) &=& c^{1/\eta} S(x,y).
 \EQNY
Furthermore, the function $ \widetilde{S}_u(1,1)$ is regularly varying at infinity with index $-1/\eta$. Other authors refer to $\eta$ as the coefficient of tail dependence (see e.g., Resnick (2002), or Reiss and Thomas (2007)).

In this paper we consider the problem of calculating the residual dependence index $\eta$ for
the bivariate random vector $(X_1,X_2)$ with stochastic representation \eqref{e} assuming that
the distribution function $F$  is in the Gumbel max-domain of attraction. We show that $\eta$ does not always exist.
In certain instances when it exists we prove that $\eta$ is defined in terms of
$\rho$ and the Weibull tail-coefficient $\theta$ (see below \eqref{eq:resid:eta}).\\
In Section 3 we propose an estimator of the residual dependence
index $\eta$. Definition, calculation
and estimation of the partial residual dependence index for
multivariate elliptical distributions are placed in Section 4.
In the multivariate setup the partial residual dependence indexes
(if they exit) are  determined by the unique solution of specific
quadratic programming problem, and the Weibull tail-coefficient
$\theta$. Proofs of all the results are relegated to Section 5 (last one).

\section{Calculation of the Residual Dependence Index}
Let $(X_1,X_2)$ be an elliptical random vector with stochastic representation \eqref{e},
and let $R$ be the positive associated random radius with distribution function $F$.
We assume in the following $F(0)=0$ and $F(x)<1, \forall x>0$. If $X_1$ and $X_2$ are standard Gaussian random variables, then it is well-known that (see e.g., Reiss and Thomas (2007))
$X_1$ and $X_2$ are asymptotically independent for any $\rho\in (-1,1)$. Furthermore, the residual dependence index $\eta$
is given by
$$ \eta:= (1+ \rho)/2\in (0,1)$$
 and (see p.\ 322 in Reiss and Thomas (2007))
\BQN\label{eq:reiss}
\widetilde{S}_u(1,1)&=& (1+o(1))  \frac{ (1- \rho^2)^{3/2}}{  (1- \rho)^2} ( 4\pi)^{- \rho/(1+\rho)}
(\ln u)^{- \rho/(1+\rho)  }  u^{ -2/(1+\rho)}, \quad u\to \IF.
\EQN
In our notation $o(1)$ means $\limit{u} o(1)=0.$

When $X_1$ and $X_2$ are independent, then $R^2\equaldis X_1^2+X_2^2$ is Chi-squared
distributed with 2 degrees of freedom. The distribution function $F$ of $R$ is in this case in the max-domain of attraction of the Gumbel distribution $\Lambda(x)=\exp(-\exp(-x)), x\inr$. From the extreme value theory we know (see e.g., Resnick (1987), Reiss (1989), Falk et al.\ (2004),
Embrechts et al.\ (1997), or de Haan and Ferreira (2006))
 that the distribution function $F$ of  $R$ is in the max-domain of attraction of $\Lambda$,
 if for some positive scaling function $w$ we have
\BQN \label{eq:rdfd}
\limit{u} \frac{1 -F(u+x/w(u))} {1- F(u)} &=&
\exp(-x),\quad \forall x\inr.
\EQN 
In the Gaussian case $\eta$ is strictly less than 1 and $w(u)=(1+o(1))u,u>0$. 
In the more general elliptical setup of this paper it turns out that interesting cases for calculation of $\eta$
are when
\BQN \label{eq:wL}
w(u)&= &u^{\theta-1} L(u),\quad \theta\in [0,\IF),
\EQN
 with $L$ a positive
slowly varying function at infinity satisfying $\limit{u}
L(cu)/L(u)=1, \forall c>0$. We refer to $\theta$ as the Weibull
tail-coefficient index (see Girard (2004)).\\ 
We show below that the elliptical model exhibits two main two main features, namely: a) the residual dependence index 
$\eta$ (when it exists) depends on both $\rho$ and $\theta$,
being in fact an increasing function of $\rho$ and $1/\theta$, and b) it is possibly to have $\eta=1$ when $\theta=0$ and 
$\limit{u} L(u)=\IF$. More interestingly, both $X,Y$ are asymptotically independent even when $\eta=1$. \\
The main result of this section is the following theorem.\\

\BT \label{theo:2} Let $(X_1,X_2),\rho \in (-1,1)$ be a bivariate
elliptical random vector with stochastic representation \eqref{e}. Assume that $R$ has distribution function $F$ which satisfies \eqref{eq:rdfd} with some positive scaling function $w$.\\
(i) Suppose that
\BQN \label{zeta2} \limit{u}
\frac{w(\alpha_{\rho} u)}{w(u)} &=& \infty, \quad \text{ with }\alpha_{\rho}:=\sqrt{2/(1+\rho)}>1
\EQN
holds, then
\BQN\label{eq:knof}
\limit{u} S_u(x,y)&=&\IF, \quad \text{ if } x>1,y>1, \quad \limit{u}
S_u(x,y)=0, \quad  \text{ if }  x,y\in (0,1).
\EQN
(ii) If for some $\theta\in [0,\IF)$ \BQN \label{zeta1}
\limit{u} \frac{w(cu)}{w(u)} &=& c^{\theta-1}, \quad \forall c>0,
\EQN
then for any $x,y\in (0,\IF)$
\BQN\label{eq:resid:eta}
\limit{u} S_u(x,y)&=&
(xy)^{1/(2\eta)} , \quad
\eta:=\fracl{1+\rho}{2}^{\theta/2}= \alpha_\rho ^{- \theta} \in (0,1],
\EQN
and $\widetilde{S}_{u}(1,1)$ is regularly varying at infinity with index $-1/\eta$.\\
(iii) Let $Q^{-1}$ denote the inverse of the distribution function of $X_1$.
As $u \to \IF$ we have the asymptotic expansion
\BQN
\widetilde{S}_u(1,1)&=& (1+o(1))\frac{\alpha_{\rho}^2 (1-\rho^2)^{3/2}}{ 2 \pi (1 - \rho )^2}
 \frac{1- F(b_*(u))}{b_*(u) w(b_*(u))}, \quad
 b_*(u):= \alpha_\rho Q^{-1}(1- 1/u).
\EQN
\ET

\begin{remarks} 
\COM{
1) Statement i) above is proved only for $x,y$ strictly larger or smaller than 1.
In fact the convergence to $\IF$ or to $0$ can be shown for the case
$xy>1$ or $xy<1$, respectively,  (where $x,y\in (0,\IF))$ by
imposing a further asymptotic condition $(u\to \IF)$ on
$w(\alpha_{\rho} u+ z/w(u))/w(\alpha_\rho u),z\inr$.
}
\COM{We note that
locally uniformly on $x\inr $ (see e.g., Resnick (1987))
$$  \frac{w(u+ x/w(u))}{w(u)} \to 1, \quad u\to \IF.$$
Furthermore, the scaling function $w$ can be defined asymptotically
by
\ass{2)} The residual dependence index $\eta$ in
\eqref{eq:resid:eta} is an increasing function of $\rho$ and
$1/\theta$.
}

1) The scaling function $w$ in \eqref{eq:rdfd} can be defined asymptotically by (see e.g., Resnick (1987))
\BQN \label{eq:w}
 w(u)&:= & \frac{(1+o(1))[1- F(u)]}{\int_{u}^{\IF} [1- F(s)]\, ds}, \quad u\to \IF.
\EQN
Further, we have
\BQN
\limit{u} u w(u)&=& \IF.
\EQN
Hence in the model \eqref{eq:wL} the Weibull tail-coefficient $\theta$ is necessarily non-negative, and if $\theta=0$, then we need to suppose further that $\limit{u} L(u)=0$. Two interesting distributions with $\theta=0$ and $L(u)= c \ln u, c\in (0,\IF), u>0$ are
Benktakder type I  and Lognormal one, see Embrechts et al.\ (1997), pp. 149-150.\\
The interesting model $\theta=0$ is kindly suggested by the referee of the paper. A striking feature of this model is that
residual dependence index $\eta$ equals 1, thus  not depending on $\rho$ at all. Furthermore, in view of Hashorva (2005a) $X$ and $Y$ are asymptotically independent.

2) In view of \eqref{eq:berm12} $X_1$ is a product of $R$ and an independent random variable $U_1$.
 Asymptotics of random products
are investigated by several authors, see for recent results Tang and Tsitsiashvili (2003,2004), Tang (2006,2008).
\end{remarks}

Next, we present three  examples.

{\bf Example 1.} Let $(X_1,X_2),\rho$ be as in \eqref{e} with associated random radius $R\sim \Lambda$. Clearly, the unit Gumbel distribution
$\Lambda$ is in the Gumbel max-domain of attraction. An admissible choice for the scaling function is $w(u)=1,\forall u>0$.
Consequently, \eqref{zeta1} holds with $\theta=1$, implying that $\widetilde{S}_u(1,1)$ is regularly varying with index
$-(0.5+ \rho/2)^{-1/2}$.\\

{\bf Example 2.} Under the setup of Example 1 we assume further that $R$ has distribution function $F$
in the Gumbel max-domain of attraction with the scaling function $w(u)= \exp(a u),u>0,$ with $a$ some positive
constant. Such $F$ exists and can easily be constructed if we assume that $F$ possesses a density function $f$,
requiring further  $f(u)/[1- F(u)]= w(u), \forall u>0$. For this choice of the scaling function $w$ \eqref{zeta2} holds.
Hence $\limit{u} S_u(x,y)=0$ for any $x,y \in (0, 1)$.\\

{\bf Example 3.} [Kotz Type III]  Again with the setup of Example 1 if for all large
$u$
\BQN\label{eq:R:K}
\pk{R> u} &=&  (1+o(1))K  u^{N}\exp(-r u^\theta), \quad K>0,\theta>0, N\inr,
 \EQN
then we refer to  $(X_1,X_2)$ as a Kotz Type III elliptical random vector.
Since we assume that $\theta$ is positive, $R$ has distribution function $F$ in the Gumbel max-domain of attraction with the scaling function
$$w(u)= (1+o(1))r \theta u ^{\theta-1}, \quad u\to \IF.$$
Consequently, \eqref{zeta1} holds and $\eta= \alpha_\rho^{-\theta}\in (0,1].$
Next, if we define  $b(u):=  Q^{-1}(1- 1/u), u>0$ with $Q$ the distribution function of $X_1$, then
\netheo{theo:2} implies
\BQNY
\widetilde{S}_{u}(1,1)&=& 
(1+o(1)) \frac{K\alpha_{\rho}^{N- \theta+2} (1-\rho^2)^{3/2}}{ 2 \pi r \theta (1 - \rho )^2}
(b(u))^{ N- \theta} \exp(- r ( \alpha_\rho b(u))^\theta ), \quad u\to \IF.
\EQNY
In view of Theorem 12.3.1 in Berman (1992) and the fact that $Q$ is symmetric about 0 (recall (1.1))
\BQNY
1- Q( u) &=& \frac{1}{2}\pk{ X_1^2> u}\\
&=& \frac{1}{2} \pk{R^2 U_1 > u} = (1+o(1)) \frac{K}{\sqrt{2 \pi r \theta    }} u ^{N-\theta/2} \exp(-ru^\theta), \quad u\to \IF,
\EQNY
where  $U_1^2$ is beta distributed with parameters $1/2,1/2$ being independent of $R$. Hence we may define $b(u)$ asymptotically as (see Embrechts et al.\ (1997))
$$ b(u)= (r^{-1}\ln u )^{1/\theta}\Biggl[ 1+ \frac{(1+o(1))}
{\theta  \ln u}\Bigl[ (N- \theta/2) \ln (r^{-1} \ln u )/\theta+ \ln K - \frac{1}{2} \ln (2 \pi r \theta) \Bigr]\Biggr]
, \quad u \to \IF.
$$
Thus we arrive at (set $\lambda:=\alpha_\rho^{\theta})$
\BQNY
\widetilde{S}_{u}(1,1)
 &=& (1+o(1))\frac{\alpha_{\rho}^{N- \theta+2} (1-\rho^2)^{3/2} r^{( \lambda-1)N/\theta } }{  K(1 - \rho )^2 }
\fracl{K^2}{2 \pi  \theta    } ^{1-\lambda/2} (\ln u)
^{ (1- \lambda)N/\theta  +  \lambda/2 -1 } u^{ -\lambda} ,  \quad u\to \IF.
\EQNY
In the special case
$$K=1, \quad r=1/2, \quad \theta=2, \quad N=0, \quad \alpha_\rho= \sqrt{2/(1+\rho)},$$
which holds in particular if both $X_1$ and $X_2$ are standard Gaussian random variables we retrieve \eqref{eq:reiss}.

\section{Estimation of $\eta$ in the Weibull Model}
In view of \netheo{theo:2} if the scaling function $w$ is  regularly varying with index $\theta-1$,
then the residual dependence index $\eta$ is defined in terms of $\rho$ and $\theta$.
Let $(X_{k1}, X_{k2}), k=1 \ldot n$ be a sample of bivariate elliptical random vectors
with stochastic representation \eqref{e} (where $\rho\in (-1,1)$ is assumed).
Then a non-parametric estimator  $\hat \rho_n$ of $\rho$ is given by (see e.g. Peng (2008), Li  and Peng (2009))
\BQN\label{eq:Bl}
 \hat \rho_n&:= &\sin(\pi \hat \tau_n /2), \quad n>1,
 \EQN
where $\hat \tau_n$ is the empirical estimator of the Kendall's tau.

Good performing estimators of the so-called Weibull tail-coefficient are the Girard and Zipf estimators, see e.g., Girard (2004).
Referring to the aforementioned paper we say that the random radius $R\sim F$ possesses  a Weibull-tail
distribution if
\BQN\label{eq:WT}
1- F(x)&=& \exp(- H(x)), \quad H^{-1}(x)=\inf\{ t: H(t) \ge x\}=  x^{1/\theta} L_1(x)
 \EQN
holds with $L_1$ a positive slowly varying function at infinity.
 Gardes and Girard (2006) and Diebolt et al.\ (2008) give several examples of Weibull-tail distributions. Prominent instances are the Gaussian, Gamma and extended Weibull distributions.\\
By the properties of slowly varying functions (see e.g., de Haan and Ferreira (2006)) we may write  \eqref{eq:WT}  alternatively as
\BQN\label{eq:WT2}
1- F(x)&=& \exp(- x^{\theta} L_2(x)),
 \EQN
where $L_2$ is another slowly varying function which is asymptotically unique.

From the estimation point of view, a tractable class of the Weibull-tail distributions is constructed when $F$
 is in the Gumbel max-domain of attraction with the scaling function $w$ defined asymptotically by
\BQN\label{eq:Abd}
w(u)& =&\frac{c u^{\theta-1}}{1+ t_1(u)}, \quad c\in (0,\IF),
\EQN
where  $t_1$ is a regularly varying function at infinity with index $\kappa:=\theta \mu, \mu \in (-\IF,  0)$ implying
\BQN \label{eq:KotzM}
1- F(u) &=& \exp(- c u^\theta(1+ t_2(u))), \quad u>0,
\EQN
where $t_2$ is another regularly varying function at infinity with index $\kappa.$

Under the assumption \eqref{eq:Abd} it follows that (see Berman (1992) or Hashorva (2005a))
$$ \pk{X_1> u}= \exp(- c u^\theta(1+ t_3(u))), \quad  u\inr,$$
with $t_3$ again a regularly varying function at infinity with index $\kappa.$

Next, assume that the associated random radius $R$ defining the random
sample $(X_{k1},X_{k2}), k=1 \ldot n,n>1$ possesses a Weibull-tail
distribution $F$ such that \eqref{eq:Abd} holds. Write $X_{1:n}\le
\cdots \le X_{n:n}$ for the associated order statistics of $X_{11}
\ldot X_{n1}$. Following Gardes and Girard (2006) we might estimate $\theta$ by
$$ \hat \theta_n:=
\frac{1}{T_n}\frac{1}{k_n}\sum_{i=1}^n \Bigl( \log X_{n-i+1:n}- \log
X_{n- k_n+1:n}\Bigr), $$
with $1 \le k_n \le n, T_n>0, n\ge 1$ given
constants satisfying
$$ \limit{n} k_n=\IF, \quad \limit{n} \frac{k_n}{n}=0, \quad
\limit{n} \log(T_n/k_n)= 1, \quad \limit{n} \sqrt{k_n}
b(\log(n/k_n))\to \lambda\inr,$$ where the function $b$ is a regularly varying function with index $\theta$ appearing in a second order
asymptotic condition imposed on $F$ (being thus related to $L_1$).

Asymptotic properties of $\hat \theta_n$ are discussed in Gardes and Girard (2006) and Diebolt et al.\ (2008).
Next, based on our main result we propose an estimator for the residual
dependence index $\eta$ given by
\BQN
\hat \eta_n&:= &\Bigl( (1+ \hat
\rho_n )/2\Bigr )^{\hat \theta_n/2},  \quad n>1.
\EQN
Asymptotic properties of $\ass{\hat \eta_n}$ follow by utilising the
asymptotic properties of both $\hat \rho_n$ and  $\hat \theta_n$.

We note in passing that the constant $c$ can be estimated by
 \BQN \hat c_{n}&=& \frac{1}{k_n} \sum_{i=1}^{k_n}
\frac{ \log (n/i) }{X_{n-i+1:n}}, \quad  n>1.
\EQN

\section{Partial Residual Dependence Index }
Consider $\X:= (X_1\ldot X_k)^\top ,k\ge 2$ an elliptical random vector in $\R^k$ with stochastic representation
\BQN \label{eq:multie}
\X&\equaldis &R A^\top \U ,
 \EQN
where $R$ is again the positive associated random radius of $\X$ with distribution function $F$ independent of
$\U:=(U_1 \ldot U_k)^\top $ which is uniformly distributed on the unit sphere of $\R^k$ and $A \inr^{k\times k}$ is
a non-singular matrix (here $^\top$ stands for the transpose sign).
The distribution function of the random vector $\X  $ is determined by the positive definite matrix
$\SI:=A^\top A$, the distribution function $F$ and the distribution function of $\vk{U}$ (which is known). See Cambanis et al.\ (1981) or Fang et al.\ (1990) for more details on elliptical
distributions.

We assume in the following
that $F(0)=0, F(x)< 1, \forall x>0$ and $\SI$ is a correlation matrix i.e.,
all the entries of the main diagonal equal 1.
If the distribution function $F$ is in the Gumbel max-domain of attraction, then each pair $X_i,X_j, i\not =j, i,j\le k$ is
asymptotically independent. If further the scaling function $w$ satisfies \eqref{zeta1},
then by Lemma 12.1.2 in Berman (1992), Proposition 3.4 in Hashorva (2005a) and \netheo{theo:2} it follows that the residual dependence index
$\eta_{ij}$ (related to $(X_i,X_j)$) is
$$ \eta_{ij}= \alpha_{\rho_{ij}}^{- \theta},$$
with $\rho_{ij}\in (-1,1)$ the $ij$-th entry of $\SI$.\\
Let $I$ be a given non-empty index subset of $\{1 \ldot k\}$. By the assumption on $\SI$ we have $X_i\equaldis X_1, i=2 \ldot k$. Next, define
$\widetilde{S}_{u,I}(\x)$ by
$$ \widetilde{S}_{u,I}(\x):= \pb{Q(X_i)> 1- \frac{x_i}{u}, \forall i\in I}, \quad u>0, \x:=(x_1 \ldot x_k)^\top \in (0,\IF)^k,$$
with $Q$ the distribution function of $X_1$.\\
If $\widetilde{S}_{u,I}(\vk{1})$ is regularly varying with index $-1/\eta_I, \eta_I \in (0, 1]$, then we refer to $\eta_I$ as the partial residual dependence  index of the subvector $\X_I:=(X_i,i\in I)^\top$, or shortly as the partial residual dependence index.

The submatrix of $\SI$ obtained by retaining the rows and the columns of $\SI$ with indices in $J$ and $I$, respectively,
(assume $I$ has less than $k$ elements) is denoted by $\SI_{JI},J:=  \{1 \ldot k\} \setminus I$. We define  similarly $\SI_{II}$ and
$\x_I$ for  $\x\inr^k$.\\
 Since in our model $\SI=A^\top A$ is positive definite the inverse matrix $\SI^{-1}_{II}$ of $\SI_{II}$
exits. Next, we write $ \alpha_{I}$ for the unique solution of the quadratic programming problem
\BQN \label{eq:min}
\text{  minimise  the objective function  } \quad  \y^\top_I \SI_{II}^{-1} \y_I , \quad  \y:= (y_1 \ldot y_{k})^\top
\in \R^{k}, \quad y_i\ge 1, \quad \forall i\in I.
\EQN
In the multivariate setup again two interesting features are observed, namely a) the partial residual index  $\eta_I$ depends on both 
$\alpha_I$ and $\theta$ (when it exists), and $b)$ when $\theta=1$, then $\eta_I=1$ for all $I\subset \njk$.  

\BT \label{theo:multi}
Let $\X$ be an elliptical random vector in $\R^k,k\ge 2,$ with stochastic representation \eqref{eq:multie}.
Assume that the associated random radius $R$ is almost surely positive with distribution function $F$
satisfying \eqref{eq:rdfd}. If the scaling function $w$ satisfies \eqref{zeta1}, then
for any non-empty index set $I\subset \{1 \ldot k\}$ with $m\le k$ elements and any $\x\in (0,\IF)^k$ we have
\BQN
\limit{u} \frac{\widetilde{S}_{u,I}(\x)}{\widetilde{S}_{u,I}(\vk{1})} &=&
\Bigl(\prod_{j\in K} x_j ^{\gamma_j} \Bigr), \quad \gamma_j:= \alpha_I^{\theta-1} (\vk{e}_j^\top \SI_{KK}^{-1} \vk{1}_K) \in (0,\IF),
\quad \forall j\in K,
\EQN
where $K$ is a unique subset of $I$ with $l>0$ elements such that $\SI^{-1}_{KK}\vk{1}_K$
is a vector with positive elements and $\vk{e}_j$ is the $j$-th unit vector in $\R^l$.
Furthermore, if $M:= I \setminus K$ is not empty, then the vector $ \SI_{KM} \SI_{MM}^{-1} \vk{1}_M - \vk{1}_K $
has non-negative components and $\widetilde{S}_{u,I}(\vk{1})$ is regularly varying with index $-1/\eta_I$ where
\BQN
\eta_I &:=& \alpha_I^{-\theta} \in \Bigl( ( \vk{1}_I^\top \SI_{II}^{-1} \vk{1}_I)^{-\theta} , 1\Bigr).
\EQN
\ET

\begin{remarks}
1) Estimation of the partial residual dependence index $\eta_I$ related to a given index set $I$
requires estimation of the attained minimum $\alpha_I$ of the related quadratic programming problem
and the Weibull tail-coefficient $\theta$. For any $i,j, i\not =j$ an
estimator of $\rho_{ij}$  (the $ij$-th entry of $\SI$) can be defined analogously to \eqref{eq:Bl} by
\BQNY
\hat \rho_{ij,n}&:= &\sin(\pi \hat \tau_{ij,n} /2), \quad n>1,
\EQNY
with $\hat \tau_{ij,n}$ the corresponding empirical estimator of the Kendall's tau.

An estimator of $\alpha_I$ can be constructed if we already have  estimated the precision matrix $\SI_{II}^{-1}$. Estimation of
$\SI^{-1}_I$ is recently discussed for Kotz distributions in Sarr and Gupta (2008).\\
If $\hat \alpha_{I,n}$ denotes an estimator of $\alpha_I$, and $\hat \theta_n$ an estimator of
the Weibull tail-coefficient, then in view of our results we can estimate the partial residual index $\eta_I$ by
\BQN
\hat \eta_{I,n}&:=&\hat \alpha_{I,n}^{-\hat \theta_n}.
\EQN
2) In the case that \eqref{zeta2} holds with $\alpha_I$ instead of $\alpha_{\rho}$, then we cannot define $\eta_I$.

3) If $\SI^{-1}_{II} \vk{1}_I$ has positive elements, then the subset $K$ in \netheo{theo:multi} equals $I$.
This is in particular the case if $I$ has only two elements, or when the non-diagonal elements of $\SI$ are all equal, say to
$\rho\in (-1,1)$.\\ If $K\not=I$, then for estimating $\alpha_I$, we also need to identify the elements of
$K$, which is not an easy task in general.

4) It is well-known that the solution of the attained minimum $\alpha_I$ of the quadratic programming problem above
is related to the exact tail asymptotics of the Gaussian random vectors, see
Dai and Mukherjea (2001),  or  Hashorva (2005b, 2007b) for more details.

5) In the particular case $w$ is given by \eqref{eq:wL} with $\theta=0$, then $\eta_I=1$ for all subsets $I$ of $\njk$.
Furthermore, $\SI$ does not influence $\eta_I$, which is in particular the case for multivariate lognormal distributions.
\end{remarks}

Next we consider the trivariate setup in some details. The following lemma gives an explicit formula
for  $\alpha_I, I=\{1,2,3\}$, which is useful for the estimation of $\alpha_I$.

\COM{
3. Letter Ref.

4. Ausdrücken.
c) Andi Quick-Winns, d) Ferien beantragen, e) Ausdrücken für morgen?
h) ticket Post, morgen, amt ..., dokumenta!!, reinigen FILES.
}

\BL \label{lem1} Let $\Sigma\inr^{3\times 3}$ be a positive definite
correlation matrix (with 1's in the main diagonal) and non-diagonal
entries $\rho_{ij}\in (-1,1), i\not=j, i,j\le 3$. Define $\rho_{min}:=\min(
\rho_{12},\rho_{13},\rho_{23})$ and set $\alpha:= \min_{x_i\ge 1,
i=1,2,3}\x^\top
\Sigma^{-1}\x$.\\
(i) If $1 + 2 \rho_{min}- \rho_{12}-\rho_{13}-\rho_{23}>0,$ then we
have (here $\vk{1}=(1,1,1)^\top)$
\BQN \label{alpha:A1}
\alpha&=& \vk{1}^\top \Sigma^{-1}\vk{1}\notag \\
&=&\frac{ 3- 2(\rho_{12}+\rho_{13}+ \rho_{23}) - \rho_{12}^2 -
\rho_{13}^2 - \rho_{23}^2
 + 2 ( \rho_{12}\rho_{13} + \rho_{12}\rho_{23}+\rho_{13}\rho_{23})}
{ 1 + 2 \rho_{12}\rho_{13} \rho_{23} -  \rho_{12}^2 - \rho_{13}^2 -
\rho_{23}^2}.
\EQN
(ii) If $1 + 2 \rho_{min}- \rho_{12}-\rho_{13}-\rho_{23}\le 0$, then there
exists a unique index set $\{i,j\} \subset \{1,2,3\}$ such that
\BQN
\rho_{min}&=&\rho_{ij}< \min_{k\not =i, k\not =l, k,l\le 3}\rho_{lk}.
\EQN
Moreover we have
\BQN\label{alpha:A2}
\alpha&=& \vk{1}^\top_K \Sigma^{-1}_{KK}\vk{1}_{K}= (1,1)^\top \Sigma^{-1}_{KK} (1,1)= \frac{2}{1+ \rho_{ij}}.
\EQN
\EL

{\bf Example 4}. [Kotz Type III, 3-dimensional  Case]. Let $\X$ be an elliptical random vector in $\R^3$
with stochastic representation \eqref{eq:multie}, where the matrix $A$ is non-singular and set $\SI:=A^\top A$. We denote by
$\rho_{ij}$ the $ij$-th entry of $\SI$. Assume that $\rho_{ii}=1, i=1 \ldot k$ and $R$ satisfies \eqref{eq:R:K}
as $u\to \IF$. Again we refer to $\X$ as a Kotz Type III random vector. In view of Lemma 12.1.2 in Berman (1992)
we have for any index set $I=\{k,l\} \subset \{1,2 3\}$ with two elements
$$ (X_k, X_l) \equaldis R\Bigl (U_1, \rho_{kl} U_1+ \sqrt{1- \rho_{kl} ^2} U_2 \Bigr),$$
where $(U_1,U_2)$ with uniform distribution on the unit
circle of $\R^2$ is independent of $R$. Hence we can estimate $\rho_{kl}$ as in \eqref{eq:Bl}. Let $\hat \rho_{12,n},\hat \rho_{13,n},\hat \rho_{23,n},n>1$ denote these
estimators.\\
Next, consider the case $I=\{1,2, 3\}$. In view of \netheo{theo:multi} and \nelem{lem1} $ \eta_I= \alpha^{-\theta},$ with $\alpha$ defined in \eqref{eq:min}.  If
$$ 1 + 2 \min (\hat \rho_{12,n},\hat \rho_{13,n},\hat \rho_{23,n})  -
\hat \rho_{12,n}- \hat \rho_{13,n}- \hat \rho_{22,n}>0,$$
then the estimator of $\alpha$ is obtained by plugging in \label{alpha:A1}
the estimators  $\hat \rho_{12,n},\hat \rho_{13,n},\hat \rho_{23,n}$. Otherwise, we estimate
$$ \hat \rho_{min,n}:= \min (\hat \rho_{12,n},\hat \rho_{13,n},\hat \rho_{23,n}   ), \quad n>1,$$
and obtain the estimator of $\alpha$ be plugging in $\hat \rho_{min,n}$ in \eqref{alpha:A2}.
The Weibull tail-coefficient $\theta$ can then be further estimated as previously discussed in Section 3.

\section{Proofs}
\def\wYu{w^*(u)}
\prooftheo{theo:2} Let $Q$ be the distribution function of $X_1$ with inverse $Q^{-1}$
($Q$ is a continuous function, see e.g., Berman (1992)). The Gumbel max-domain of attraction assumption on $F$
implies (see e.g., Reiss (1989), or de Haan and Ferreira (2006))
\BQN \label{eq:dehaan}
 w(b(u))[ Q^{-1}(1- x/u)- b(u)]  &\to & - \ln x, \quad u \to \IF
 \EQN
locally uniformly for $x\in (0,\IF)$, with $b(u):=Q^{-1}(1- 1/u), u>0$. Next set
$$ w^*(u):=w( \alpha_{\rho} b(u) ),u>0, \quad \alpha_{\rho}:= \sqrt{2/(1+\rho)}>1, \quad \rho \in (-1,1). $$
For any $u, x,y $ positive
we may further write (recall $X_1 \equaldis X_2)$
\BQNY
\widetilde{S}_u(x,y)&=& \pb{ Q(X_1)> 1- \frac{x}{u}, Q(X_2)> 1- \frac{y}{u}}\\
 &=& \pb{ X_1> Q^{-1}(1- \frac{x}{u}), X_2> Q^{-1}(1- \frac{y}{u})}\\
&=&\pb{ X_1> b(u)- (1+o(1))\frac{\ln x}{w(b(u))}, X_2> b(u)- (1+o(1)) \frac{\ln y}{w(b(u))}}.
\EQNY
In view of Theorem 5 in Hashorva (2007a) for any $s,t$ positive
$$ \limit{u} \pb{\wYu(X_1- b(u))> s, \wYu(X_2- b(u))> t
\Bigl \lvert X_1> b(u), X_2> b(u)}= \pk{X'_1>s, X'_2>t}, $$
holds with $X_1',X_2'$ two independent exponentially distributed random variables with
mean $\lambda_{\rho}:=\sqrt{2(1+ \rho)}$. Hence if $x,y\in (0,1)$, then $-\ln x,-\ln y\in (0,\IF)$,
thus \eqref{eq:knof} follows easily. For any $x>1$ and $y>1$ we may write
$$ \limit{u} S_u(x,y)= \limit{u} \frac{1}{S_u(1/x,1/y)}= \IF.$$
Next, if \eqref{zeta1} holds, then
$$  \frac{ \wYu }{\alpha_{\rho}^{\theta-1}  \ln x}\Biggl[ G^{-1}(1- x/u)- b(u) \Biggr]\to -1, \quad   u \to \IF$$
locally uniformly for any $x>0$. Consequently, with the same arguments as above
for any $x,y\in (0,1]$ we obtain
\BQNY
 \limit{u} S_u(x,y)&=& \pk{X'_1>- \alpha_{\rho}^{\theta-1} \ln x , X_2'>- \alpha_{\rho}^{\theta-1}\ln y}\\
 &=&\exp\Bigl( \frac{\alpha_{\rho}^{\theta-1} }{\lambda_{\rho}} \ln (xy)\Bigr)=
 \exp\Bigl( \frac{\alpha_{\rho}^{\theta} }{2} \ln (xy)\Bigr) =:S(x,y).
 \EQNY
The result for $x \in (1,\IF)$ and $y$ positive, or $x$ positive and $y\in (1,\IF)$ as well as the statement (iii)
can be now established using directly Theorem 2 in the aforementioned paper. Since for any $c,x,y$ positive
$$ S(cx,cy)= S(x,y) c^{   1/ \eta},$$
with
$$ \eta:= \alpha_{\rho}^{-\theta }=\fracl{1+\rho}{2}^{\theta/2}\in (0, 1] $$
the result follows.
\QED

\prooftheo{theo:multi} Let $I$ be a no-empty subset of $\{1 \ldot k\}$ with $m\le k$ elements.
The random vector $\X_I:=(X_i, i\in I)^\top$  is again an elliptical random vector with stochastic representation
(Cambanis et al.\ (1981))
$$ \X_I\equaldis R_I B \vk{V},$$
with positive associate random radius $R_I$,  square matrix $B$ such  that $B^\top B= \SI_{II}$ and
$\vk{V}$ uniformly distributed on the unit sphere of $\R^{m}$ being independent of $R_I$. As shown in
Berman (1992) the associated random radius $R_I$ has distribution function $F_I$ in the
Gumbel max-domain of attraction with the same scaling function $w$ as $F$ the distribution function of $R$.
 By Proposition 2.1 in Hashorva (2005b) there exist a unique subset $K\subset I$ with $l> 0$ elements
 such that
 $$ \alpha_I:=\min_{\y\inr^m, y_i\ge 1, i=1 \ldot m} \y^\top\SI^{-1}_{II} \y=
  \vk{1}_K^\top \SI_{KK}^{-1} \vk{1}_K>0,$$
$\SI_{KK}^{-1}\vk{1}_K$ has non-negative components, and if $M:=I\setminus K$ is not empty,
then  $\Sigma_{KM} \SI_{MM}^{-1} \vk{1}_M - \vk{1}_K$ has non-negative elements.\\
As in the proof of \netheo{theo:2} for any $\x\in (0,\IF)^k$ applying further Theorem  3.4 in Hashorva (2007b)
we obtain
\BQNY
\limit{u} \frac{ \widetilde{S}_{u,I}(\x) }{ \widetilde{S}_{u,I}(\vk{1})}
&=&
\Bigl( \prod_{j \in K}  x_j^{\mu_J  \alpha_I^{\theta-1}}\Bigr )=:S(\x)\in (0,\IF),
\EQNY
with $\mu_j:= \vk{e}_j^\top \SI^{-1}_{KK} \vk{1}_K> 0,$ and $\vk{e}_j$ the $j$-th unit vector in $\R^{l}$.
We have further
$$ \sum_{j\in K} \mu_j= \sum_{j\in K} \vk{e}_j^\top \SI^{-1}_{KK} \vk{1}_K= \alpha_I,$$
hence for any $c>0$ and any $\x=(x_1 \ldot x_k)^\top \in (0,\IF)^k$ we may write
\BQNY
\frac{S(c x_1 \ldot c x_k)}{S( x_1 \ldot  x_k)  }&= &\Bigl( \prod_{j\in K} c^{ \mu_j  \alpha_I^{\theta-1} }\Bigr)\\
& =&c^{\alpha_I^{\theta-1} \sum_{j\in K} \mu_j}= c^{\alpha_I^{\theta}}.
\EQNY
Consequently,  $\eta_I= \alpha_{I}^{- \theta}$, thus the result follows. \QED

\prooflem{lem1} The proof of the first statement is shown in
Lemma 3.2 in Hashorva and H\"usler (2002). We show next the second statement. Assume therefore that
 $$1 + 2 \rho_{min}-\rho_{12}-\rho_{13}-\rho_{23}\le 0, \quad \text{ and } \rho_{min}=\rho_{12}.$$
Since $1  - \max(\rho_{23},\rho_{13})> 0$,  then $\rho_{12}=\rho_{23}$ or $\rho_{12}=\rho_{13}$ is not possible, hence $\rho_{min}= \rho_{12}$. In view of the aforementioned lemma $\alpha= 2/(1+ \rho_{12})$, thus the result follows. \QED

{\bf Acknowledgement:} I am in debt to the referee for deep professional review, very kind report, several suggestions and corrections.

\end{document}